\documentstyle[11pt,amssymb]{article}
\begin{document}

{\Large

\noindent{\bf On $q$-orthogonal polynomials, dual to little}

\noindent{\bf and big $q$-Jacobi polynomials} }
\bigskip

\noindent{\sc N. M. Atakishiyev${}^1$ and A. U. Klimyk${}^{1,2}$}
\bigskip

\noindent ${}^1$Instituto de Matem\'aticas, UNAM, CP 62210
Cuernavaca, Morelos, M\'exico

\noindent ${}^2$Bogolyubov Institute for Theoretical Physics,
03143 Kiev, Ukraine

\medskip

E-mail: natig@matcuer.unam.mx and anatoliy@matcuer.unam.mx
\bigskip

\begin{abstract}
This paper studies properties of $q$-Jacobi polynomials and their
duals by means of operators of the discrete series representations
for the quantum algebra $U_q({\rm su}_{1,1})$. Spectra and
eigenfunctions of these operators are found explicitly. These
eigenfunctions, when normalized, form an orthonormal basis in the
representation space. The initial $U_q({\rm su}_{1,1})$-basis and
the bases of these eigenfunctions are interconnected by matrices,
whose entries are expressed in terms of little or big $q$-Jacobi
polynomials. The orthogonality by rows in these unitary connection
matrices leads to the orthogonality relations for little and big
$q$-Jacobi polynomials. The orthogonality by columns in the
connection matrices leads to an explicit form of orthogonality
relations on the countable set of points for ${}_3\phi_2$ and
${}_3\phi_1$ polynomials, which are dual to big and little
$q$-Jacobi polynomials, respectively. The orthogonality measure
for the dual little $q$-Jacobi polynomials proves to be extremal,
whereas the measure for the dual big $q$-Jacobi polynomials is not
extremal.

\end{abstract}

\bigskip

{\bf Key words.} Orthogonal $q$-polynomials, little $q$-Jacobi
polynomials, big $q$-Jacobi polynomials, Leonard pairs,
orthogonality relations, quantum algebra
\medskip

{\bf AMS subject classification.} 33D80, 33D45, 17B37
\bigskip

\noindent{\bf 1. Introduction}
\bigskip

The appearance of quantum groups and quantized universal
enveloping algebras (quantum algebras) and development of their
representations led to their applications in the theory of
$q$-orthogonal polynomials and $q$-special functions (see, for
example, [1--4]). Since the theory of quantum groups and their
representations is much more complicated than the Lie theory, the
corresponding applications are more difficult. At the first stage
of such applications, the compact quantum groups and their finite
dimensional representations have been used.

It is known that representations of the noncompact Lie group
$SU(1,1)\sim SL(2,{\Bbb R})$ are very productive for the theory of
orthogonal polynomials and special functions (see, for example,
[5], Chapter 7). Unfortunately, there are difficulties with a
satisfactory definition of the noncompact quantum group
$SU_q(1,1)$, which would give us a possibility to use such quantum
group extensively for deeper understanding the theory of
orthogonal polynomials and special functions. For this reason,
representations of the corresponding quantum algebra $U_q({\rm
su}_{1,1})$ have been commonly used for such purposes (see, for
example, [6--12]).

In this paper we continue to use representations of the positive
discrete series of the quantum algebra $U_q({\rm su}_{1,1})$ for
exploring properties of $q$-orthogonal polynomials [6--8]. In
fact, we deal with certain operators in these representations and
do not touch the Hopf structure of the algebra $U_q({\rm
su}_{1,1})$. Our study of these polynomials is related to
representation operators, which can be represented by a Jacobi
matrix. Namely, we consider those representation operators $A$,
which correspond to some particular Jacobi matrices. In the case
under discussion we diagonalize these selfadjoint bounded
operators with the aid of big or little $q$-Jacobi polynomials. An
explicit form of all eigenfunctions of these operators is found.
Since the spectra are simple, eigenfunctions of each such operator
form an orthogonal basis in the representation space. One can
normalize this basis. This normalization is effected by means of
the second operator $J$, which is related (in some sense) to
$q$-difference equations for little and big $q$-Jacobi
polynomials. As a result of the normalization, for each operator
$A$ (one of them is related to little $q$-Jacobi polynomials and
another one to big $q$-Jacobi polynomials) two orthonormal bases
in the representation space emerge: the canonical (or the initial)
basis and the basis of eigenfunctions of the operator $A$. They
are interrelated by a unitary matrix $U$ whose entries $u_{mn}$
are explicitly expressed in terms of little or big $q$-Jacobi
polynomials. Since the matrix $U$ is unitary (and in fact it is
real in our case), there are two orthogonality relations for its
elements, namely
$$
\sum _{n} u_{mn}u_{m'n}=\delta_{mm'},\ \ \ \sum _{m}
u_{mn}u_{mn'}=\delta_{nn'}.   \eqno (1.1)
$$
The first relation expresses the orthogonality relation for little
or big $q$-Jacobi polynomials. So, the orthogonality of $U$ yields
an algebraic proof of orthogonality relations for these
polynomials. In order to interpret the second relation, we
consider little and big $q$-Jacobi polynomials $P_n(q^{-m})$ as
functions of $n$. In this way one obtains two sets of orthogonal
functions (one for little $q$-Jacobi polynomials and another for
big $q$-Jacobi polynomials), which are expressed in terms of
$q$-orthogonal polynomials (which can be considered as a dual sets
of polynomials with respect to little and big $q$-Jacobi
polynomials; such duality is well known in the case of
polynomials, orthogonal on a finite set of points). The second
relation in (1.1) leads to the orthogonality relations for these
$q$-orthogonal polynomials on nonuniform lattices.

In fact, this idea extends the notion of the duality of
polynomials, orthogonal on a finite set, to the case of
polynomials, orthogonal on an infinite set of points. We have
already used this idea in [6] and [8] to show that
Al-Salam--Carlitz II polynomials are dual with respect to little
$q$-Laguerre polynomials and $q$-Meixner polynomials are dual to
big $q$-Laguerre polynomials. We emphasize at this point that
there are known theorems on dual orthogonality properties of
$q$-polynomials, whose weight functions are supported on a
discrete set of points (see, for example, [13] and [14]). However,
they are formulated in terms of orthogonal functions (see (5.1)
and (8.1) below for their explicit forms in the case of little and
big $q$-Jacobi polynomials, respectively) as dual objects with
respect to given orthogonal polynomials. Therefore, one still
needs to make one step further in order to single out an
appropriate family of dual polynomials from these functions. So,
our main motivation for this paper is to show explicitly how to
accomplish that for little and big $q$-Jacobi polynomials.

The orthogonality measure for polynomials, dual to little
$q$-Jacobi polynomials, is extremal, that is, these polynomials
form a complete set in the space $L^2$ with respect to their
orthogonality measure.

The orthogonality measure for polynomials, dual to big $q$-Jacobi
polynomials, is not extremal: these polynomials do not form a
complete set in the corresponding space $L^2$. We have found the
complementary set of orthogonal functions in this space $L^2$.
These functions are expressed in terms of the same polynomials but
with different values of parameters.

The dual little $q$-Jacobi polynomials and the dual big $q$-Jacobi
polynomials, as other well-known $q$-orthogonal polynomials, can
be applied in different branches of science. For example, they may
be useful for studying a certain type of $q$-difference equations,
which appear in applications in engineering and physics.

Throughout the sequel we always assume that $q$ is a fixed
positive number such that $q<1$. We use (without additional
explanation) notations of the theory of special functions and the
standard $q$-analysis (see, for example, [15] and [16]). In
particular, we adopt for $q$-numbers $[a]_q$ the form
$$
[a]_q := \frac{q^{a/2}-q^{-a/2}}{q^{1/2}-q^{-1/2}}, \eqno (1.2)
$$
where $a$ is any complex number and $q^{1/2}$ is assumed to be
positive. We shall also use the well-known shorthand
$(a_1,\cdots,a_k;q)_n:=(a_1;q)_n\cdots (a_k;q)_n$.
\bigskip

\noindent{\bf 2. Discrete series representations of $U_q({\rm
su}_{1,1})$}
\bigskip

The quantum algebra $U_q({\rm su}_{1,1})$ is defined as the
associative algebra, generated by the elements $J_+$, $J_-$,
$q^{J_0}$ and $q^{-J_0}$, subject to the commutation relations
$$
q^{J_0}q^{-J_0}=q^{-J_0}q^{J_0}=1,\ \ \ \ q^{J_0}J_{\pm}\,
q^{-J_0}=q^{\pm 1}J_{\pm},\ \ \ \ [J_-,J_+] =
{q^{J_0}-q^{-J_0}\over q^{1/2}-q^{-1/2}}  ,
$$
and the involution relations $(q^{J_0})^* = q^{J_0}$ and $J_+^* =
J_-$. (Observe that here we have replaced $J_-$ by $-J_-$ in the
common definition of the algebra $U_q({\rm sl}_{2})$.) For
brevity, in what follows we denote the algebra $U_q({\rm
su}_{1,1})$ by ${\rm su}_q(1,1)$.

If the algebra ${\rm su}_q(1,1)$ is realized in terms of the
operators, one may consider also the operator $J_0$. In this case
instead of the commutation relations between $J_+$, $J_-$,
$q^{J_0}$ and $q^{-J_0}$, it is convenient to work with more
familiar relations
$$
[J_0,J_\pm ]=\pm \, J_\pm ,\ \ \ \ [J_-,J_+] =
{q^{J_0}-q^{-J_0}\over q^{1/2}-q^{-1/2}} .
$$
Then the involution relations reduce to the following ones:
$$
J_0^* = J_0,\ \ \ \ J_+^* = J_-. \eqno (2.1)
$$

We are interested in the discrete series representations of ${\rm
su}_q(1,1)$ with lowest weights. These irreducible representations
will be denoted by $T^+_l$, where $l$ is a lowest weight, which
can be any positive number (see, for example, [17]).

The representation $T^+_{l}$ can be realized on the space ${\cal
L}_{l}$ of all polynomials in $x$. We choose a basis for this
space, consisting of the monomials
$$
f^{l}_n(x) := c^{l}_n\, x^n, \ \ \  n = 0,1,2,\cdots , \eqno(2.2)
$$
where
$$ c^{l}_0 = 1, \qquad c^{l}_n = \prod_{k=1}^n\,{[2l+k-1]_q^{1/2}
\over [k]_q^{1/2}} = q^{(1-2l)n/4}{(q^{2l};q)_n^{1/2}\over
(q;q)_n^{1/2}}\,, \ \ n = 1,2,3, \cdots,\eqno(2.3)
$$
and $(a;q)_n=(1-a)(1-aq)\ldots (1-aq^{n-1})$. The representation
$T^+_{l}$ is then realized by the operators
$$
J_0 =x{d\over dx} + l,\qquad J_{\pm} = x^{\pm 1}[ J_0(x)\pm l]_q
\,.
$$
As a result of this realization, we have
$$
J_+\,f^{l}_n =\sqrt{[2l+n]_q \,[n+1]_q }\, f^{l}_{n+1}
=\frac{q^{-(n+l-1/2)/2}}{1-q} \sqrt{(1-q^{n+1})(1-q^{2l+n})}
f^{l}_{n+1}, \eqno (2.4)
$$ $$
J_-\, f^{l}_n =\sqrt{[2l+n-1]_q \,[n]_q }\,f^{l}_{n-1} =
\frac{q^{-(n+l-3/2)/2}}{1-q} \sqrt{(1-q^{n})(1-q^{2l+n-1})}
f^{l}_{n-1},  \eqno (2.5)
$$  $$
J_0\, f^{l}_n = (l + n)\,f^{l}_n . \eqno (2.6)
$$

We know that the discrete series representations $T_l$ can be
realized on a Hilbert space, on which the adjointness relations
(2.1) are satisfied. In order to obtain such a Hilbert space, we
assume that the monomials $f^{l}_n(x)$, $n=0,1,2,\cdots$,
constitute an orthonormal basis for this Hilbert space. This
introduces a scalar product $\langle \cdot , \cdot \rangle$ into
the space ${\cal L}_l$. Then we close this space with respect to
this scalar product and obtain the Hilbert space, which will be
denoted by ${\cal H}_l$.

The Hilbert space ${\cal H}_l$ consists of functions (series)
$$
f(x)=\sum _{n=0}^\infty b_nf^l _n(x)=\sum _{n=0}^\infty b_nc^l_n
x^n= \sum _{n=0}^\infty a_nx^n,
$$
where $a_n=b_nc^l_n$. Since $\langle f^l_m ,f^l_n\rangle
=\delta_{mn}$ by definition, for $f(x)=\sum _{n=0}^\infty a_nx^n$
and $f'(x)=\sum _{n=0}^\infty a'_n\,x^n$ we have $\langle f ,
f'\rangle = \sum _{n=0}^\infty a_n \, \overline{ a'_n}/|c^l_n|^2$,
that is, the Hilbert space ${\cal H}_l$ consists of analytical
functions $f(x)=\sum _{n=0}^\infty a_n\,x^n$, such that $ \Vert
f\Vert ^2 \equiv \sum _{n=0}^\infty |a_n/c_n^l|^2 < \infty$.
\bigskip

\noindent{\bf 3. Leonard pair $(I_1,J)$}
\bigskip

Let $I_1$ be the operator
$$
I_1 := -\alpha q^{J_0/4}(J_+\,A +A\,  J_-) q^{J_0/4} + B \eqno
(3.1)
$$
of the representation $T^+_l$, where $\alpha =
(a/q)^{1/2}\,(1-q)$, $a=q^{2l-1}$, and $A$ and $B$ are operators
of the form
$$
A=\frac{q^{J_0-l+1/2} \sqrt{(1-bq^{J_0-l+1})(1-abq^{J_0-l+1})}}
{(1-abq^{2J_0-2l+2})\sqrt{(1-abq^{2J_0-2l+1})(1-abq^{2J_0-2l+3})}}\,,
$$   $$
B=\frac{q^{J_0-l}}{1-abq^{2J_0-2l+1}} \left(
\frac{(1-aq^{J_0-l+1})(1-abq^{J_0-l+1})}{1-abq^{2J_0-2l+2}}+a
 \frac{(1-q^{J_0-l})(1-bq^{J_0-l})}{1-abq^{2J_0-2l}} \right) .
$$
Since the bounded operator $q^{J_0}$ is diagonal in the basis $\{
f^l_n\}$, the operators $A$ and $B$ are well defined.

We have the following expression for the action of the operator
$I_1$ in the canonical basis $f^l_n$, $n=0,1,2,\cdots$:
$$
I_1\, f^l_n=-a_nf^l_{n+1}-a_{n-1}f^l_{n-1}+b_nf^l_n , \eqno (3.2)
$$
where
$$
a_n=a^{1/2}q^{n+1/2}\frac{ \sqrt{(1-q^{n+1})(1-aq^{n+1})
(1-bq^{n+1})(1-abq^{n+1})}}{(1-abq^{2n+2})\sqrt{(1-abq^{2n+1})
(1-abq^{2n+3})}} ,
$$  $$
b_n=\frac{q^n}{1-abq^{2n+1}}\left(
\frac{(1-aq^{n+1})(1-abq^{n+1})}{1-abq^{2n+2}} +a
\frac{(1-q^{n})(1-bq^{n})}{1-abq^{2n}} \right) .
$$
In order to assure that expressions for $a_n$ and $b_n$ are well
defined, we suppose that $b<q^{-1}$. Note that
$0<a=q^{2l-1}<q^{-1}$ and $l$ takes any positive value. Since the
$q^{\pm J_0}$ are symmetric operators, the operator $I_1$ is also
symmetric.

Since $a_n\to 0$ and $b_n\to 0$ when $n\to \infty$, the operator
$I_1$ is bounded. Therefore, we assume that it is defined on the
whole representation space ${\cal H}_l$. For this reason, $I_1$ is
a selfadjoint operator. Let us show that $I_1$ is a trace class
operator (we remind that a bounded selfadjoint operator is a trace
class operator if a sum of its matrix elements in an orthonormal
basis is finite; a spectrum of such an operator is discrete, with
a single accumulation point at $0$ ). For the coefficients $a_n$
and $b_n$ from (3.2), we have
$$
a_{n+1}/a_n \to q,\ \ b_{n+1}/b_n \to q \ \ \ {\rm when}\ \ \ n\to
\infty .
$$
Therefore, for the sum of all matrix elements of the operator
$I_1$ in the canonical basis we have $\sum _n (2a_n+b_n)< \infty$.
This means that $I_1$ is a trace class operator. Thus, a spectrum
of $I_1$ is discrete and have a single accumulation point at 0.
Moreover, a spectrum of $I_1$ is simple, since $I_1$ is
representable by a Jacobi matrix with $a_n\ne 0$ (see [18],
Chapter VII).

To find eigenfunctions $\xi_\lambda (x)$ of the operator $I_1$,
$I_1 \xi_\lambda (x)=\lambda \xi_\lambda (x)$, we set
$$
\xi_\lambda (x)=\sum _{n=0}^\infty \beta_n(\lambda)f^l_n (x).
$$
Acting by the operator $I_1$ upon both sides of this relation, one
derives that
$$
\sum _{n=0}^{\infty}\, \beta_n(\lambda)\,
(a_nf^l_{n+1}+a_{n-1}f^l_{n-1}-b_nf^l_n) = - \lambda
\sum_{n=0}^{\infty}\, \beta_n(\lambda)f^l_n ,
$$
where $a_n$ and $b_n$ are the same as in (3.2). Collecting in this
identity all factors, which multiply $f^l_n$ with fixed $n$, one
derives the recurrence relation for the coefficients
$\beta_n(\lambda)$:
$$
\beta_{n+1}(\lambda)a_n +\beta_{n-1}(\lambda)a_{n-1}-
\beta_{n}(\lambda)b_n= -\lambda \beta_{n}(\lambda).
$$
Making the substitution
$$
\beta_{n}(\lambda)=\left( \frac{(abq,aq;q)_n\,(1-abq^{2n+1})}
{(bq,q;q)_n\, (1-abq)(aq)^n}\right) ^{1/2} \beta'_{n}(\lambda)
$$
reduces this relation to the following one
$$
A_n \beta'_{n+1}(\lambda)+ C_n \beta'_{n-1}(\lambda)
-(A_n+C_n)\beta'_{n}(\lambda)=-\lambda \beta'_{n}(\lambda)
$$
with
$$
A_n=\frac{q^n(1-aq^{n+1})(1-abq^{n+1})}{(1-abq^{2n+1})
(1-abq^{2n+2})}, \ \ \ \
C_n=\frac{aq^n(1-q^{n})(1-bq^{n})}{(1-abq^{2n}) (1-abq^{2n+1})}.
$$
This is the recurrence relation for the little $q$-Jacobi
polynomials
$$
p_n(\lambda ;a,b|q):={}_2\phi_1 (q^{-n}, abq^{n+1};\; aq; \;
q,q\lambda )   \eqno (3.3)
$$
(see, for example, formula (7.3.1) in [15]). Therefore,
$\beta'_n(\lambda )=p_n(\lambda ;a,b|q)$ and
$$
\beta_n(\lambda )= \left( \frac{(abq,aq;q)_n\,(1-abq^{2n+1})}
{(bq,q;q)_n\, (1-abq)(aq)^n}\right)^{1/2}p_n(\lambda ;a,b|q).
\eqno (3.4)
$$
For the eigenfunctions $\xi _\lambda(x)$ we have the expression
$$
\xi _\lambda(x)=\sum_{n=0}^\infty\,\left( \frac{(abq,aq;q)_n\,
(1-abq^{2n+1})} {(bq,q;q)_n\,
(1-abq)(aq)^n}\right)^{1/2}p_n(\lambda ;a,b|q)\,f^l_n(x)
$$   $$
=\sum_{n=0}^\infty\,a^{-n/4}\,\frac{(aq;q)_n}{(q;q)_n}\,
\,\left(\frac{(abq;q)_n\,(1-abq^{2n+1})} {(bq;q)_n\,
(1-abq)(aq)^n}\right)^{1/2} p_n(\lambda ;a,b|q) x^n. \eqno (3.5)
$$
Since the spectrum of the operator $I_1$ is discrete, only a
discrete set of these functions belongs to the Hilbert space
${\cal H}_l$. This discrete set of functions determines a spectrum
of $I_1$.

Now we look for a spectrum of the operator $I_1$ and for a set of
polynomials, dual to the little $q$-Jacobi polynomials. To this
end we use the action of the operator
$$
J:= q^{-J_0+l} + ab\,q^{J_0-l+1}
$$
upon the eigenfunctions $\xi _\lambda(x)$, which belong to the
Hilbert space ${\cal H}_l$. In order to find how this operator
acts upon these functions, one can use the $q$-difference equation
$$
(q^{-n}+abq^{n+1})\,p_n(\lambda)= a\lambda^{-1}(bq\lambda
-1)\,p_{n} (q\lambda
)+\lambda^{-1}(1+a)\,p_n(\lambda)+\lambda^{-1}
(\lambda-1)\,p_n(q^{-1}\lambda) \eqno(3.6)
$$
for the little $q$-Jacobi polynomials $p_n(\lambda)\equiv
p_n(\lambda ;a,b|q)$ (see, for example, formula (3.12.5) in [19]).
Multiply both sides of (3.6) by $d_n\,f^l_n(x)$ and sum up over
$n$, where $d_n$ are the coefficients of $p_n(\lambda ;a,b|q)$ in
the expression (3.4) for the coefficients $\beta_n(\lambda)$.
Taking into account the first part of formula (3.5) and the fact
that $Jf^l_n(x)=(q^{-n}+abq^{n+1})f^l_n(x)$, one obtains the
relation
$$
J\,\xi _{\lambda}(x)= a\lambda^{-1}(bq\lambda -1)\,\xi
_{q\lambda}(x) + \lambda^{-1}(1+a)\, \xi _{\lambda}(x)+
\lambda^{-1}(\lambda-1)\, \xi_{q^{-1}\lambda}(x). \eqno (3.7)
$$

It will be shown in the next section that the spectrum of the
operator $I_1$ consists of the points $\lambda=q^n$,
$n=0,1,2,\cdots$. Thus, we see that the pair of the operators
$I_1$ and $J$ form a Leonard pair (see [20], where P.~Terwilliger
has actually introduced this notion in an effort to interpret the
results of D. Leonard [21]). Recall that a pair of operators $R_1$
and $R_2$, acting in a linear space ${\cal L}$, is a Leonard pair
if

(a) there exists a basis in ${\cal L}$, with respect to which the
operator $R_1$ is diagonal, and the operator $R_2$ has the form of
a Jacobi matrix;

(b) there exists another basis of ${\cal L}$, with respect to
which the operator $R_2$ is diagonal, and the operator $R_1$ has
the form of a Jacobi matrix.

Properties of Leonard pairs of operators in finite dimensional
spaces are studied in detail (see, for example, [22] and [23]).
Leonard pairs in infinite dimensional spaces are more complicated
and only some isolated results are known in this case (see, for
example, [23]).
\bigskip

\noindent{\bf 4. Spectrum of $I_1$ and orthogonality of little
$q$-Jacobi polynomials}
\medskip

The aim of this section is to find, by using the Leonard pair
$(I_1,J)$, a basis in the Hilbert space ${\cal H}_l$, which
consists of eigenfunctions of the operator $I_1$ in a normalized
form, and to derive explicitly the unitary matrix $U$, connecting
this basis with the canonical basis $f^l_n$, $n=0,1,2,\cdots$, in
${\cal H}_l$. This matrix directly leads the orthogonality
relation for the little $q$-Jacobi polynomials.

A word of explanation is in order. It is a well-known fact that
the little $q$-Jacobi polynomials (3.3) represent a particular
case of the big $q$-Jacobi polynomials (see formula (6.3) below)
with the vanishing parameter $c$. Therefore, it may seem that one
should start directly with the latter set of polynomials --
wherein the corresponding results for the former family (3.3) are
then recovered by assuming that $c=0$. However, we are interested
in finding dual families with respect to polynomials (3.3) and
(6.3), and establishing their properties. It turns out that an
interrelation between dual little $q$-Jacobi and dual big
$q$-Jacobi polynomials is more intricate than that occurring
between polynomials (3.3) and (6.3) themselves (for instance, the
orthogonality measure for dual little $q$-Jacobi polynomials
proves to be an extremal one, whereas the orthogonality measure
for dual big $q$-Jacobi polynomials is not extremal). For this
reason we find it more instructive to begin our discussion with
the simpler case (3.3) in this section and then to consider duals
of the polynomials (6.3) in sections 6--8.

Let us analyze a form of the spectrum of the operator $I_1$ from
the point of view of the representations of the algebra ${\rm
su}_q(1,1)$. If $\lambda$ is a spectral point of the operator
$I_1$, then (as it is easy to see from (3.7)) a successive action
by the operator $J$ upon the function (eigenfunction of $I_1$)
$\xi_\lambda$ leads to the functions
$$
\xi_{q^m\lambda}, \ \ \ m=0,\pm 1, \pm 2,\cdots . \eqno (4.1)
$$
However, since $I_1$ is a trace class operator, not all these
points can belong to the spectrum of $I_1$, since $q^{-m}\lambda
\to\infty$ when $m\to \infty$ if $\lambda\ne 0$. This means that
the coefficient $\lambda' -1$ of $\xi _{q^{-1}\lambda'}(x)$ in
(3.7) must vanish for some eigenvalue $\lambda'$. Clearly, it
vanishes when $\lambda' =1$. Moreover, this is the only
possibility for vanishing a coefficient at $\xi
_{q^{-1}\lambda'}(x)$ in (3.7), that is, the point $\lambda =1$ is
a spectral point for the operator $I_1$. Let us show that the
corresponding eigenfunction $\xi _{1}(x)\equiv \xi_{q^{0}}(x)$
belongs to the representation space ${\cal H}_l$.

Observe that by formula (II.6) of Appendix II in [15], one has
$$
p_n(1 ;a,b|q)={}_2\phi_1 (q^{-n}, abq^{n+1};\; aq; \; q,q
)=\frac{(b^{-1}q^{-n};q)_n}{(aq;q)_n} (abq^{n+1})^n.
$$
Since $(b^{-1}q^{-n};q)_n=(bq;q)_n
(-b^{-1}q^{-1})^nq^{-n(n-1)/2}$, this means that
$$
p_n(1 ;a,b|q)=\frac{(bq;q)_n}{(aq;q)_n}(-a)^nq^{n(n+1)/2}.
$$
Therefore,
$$
\langle \xi_1(x),\xi_1(x)\rangle = \sum_{n=0}^\infty
\frac{(abq,aq;q)_n\,(1-abq^{2n+1})}{(bq,q;q)_n\,
(1-abq)(aq)^n}\,p^2_n(1 ;a,b|q)
$$   $$
= \sum_{n=0}^\infty
\frac{(abq,bq;q)_n\,(1-abq^{2n+1})}{(aq,q;q)_n(1-abq)}\,
a^nq^{n^2} = \frac{(abq^2;q)_\infty}{(aq;q)_\infty} . \eqno (4.2)
$$
The last leg of this equality is obtained from formula (A.1) of
Appendix. Thus, the series (4.2) converges and, therefore, the
point $\lambda =1$ actually belongs to the spectrum of the
operator $I_1$.

Let us find other spectral points of the operator $I_1$ (recall
that a spectrum of $I_1$ is discrete). Setting $\lambda = 1$ in
(3.7), we see that the operator $J$ transforms $\xi _{q^0}(x)$
into a linear combination of the functions $\xi _{q}(x)$ and
$\xi_{q^0}(x)$. Moreover, $\xi_q(x)$ belongs to the Hilbert space
${\cal H}_l$, since the series
$$
\langle \xi _{q} ,\xi _{q} \rangle = \sum_{n=0}^\infty
\frac{(abq,aq;q)_n\,(1-abq^{2n+1})}{(bq,q;q)_n\, (1-abq)(aq)^n}
p^2_n(q ;a,b|q) <\infty
$$
is majorized by the corresponding series for $\xi_{q^0}(x)$,
considered above. Therefore, $\xi _{q}(x)$ belongs to the Hilbert
space ${\cal H}_l$ and the point $q$ is an eigenvalue of the
operator $I_1$. Similarly, setting $\lambda=q$ in (3.7), we find
that $\xi _{q^2}(x)$ is an eigenfunction of $I_1$ and the point
$q^2$ belongs to the spectrum of $I_1$. Repeating this procedure,
we find that $\xi _{q^n}(x)$, $n=0,1,2,\cdots$, are eigenfunctions
of $I_1$ and the set $q^n$, $n=0,1,2,\cdots$, belongs to the
spectrum of $I_1$. So far, we do not know yet whether other
spectral points exist or not.

The functions $\xi _{q^n}(x)$, $n=0,1,2,\cdots$, are linearly
independent elements of the representation space ${\cal H}_l$
(since they correspond to different eigenvalues of the selfadjoint
operator $I_1$). Suppose that values $q^n$, $n=0,1,2,\cdots$,
constitute a whole spectrum of the operator $I_1$. Then the set of
functions $\xi _{q^n}(x)$, $n=0,1,2,\cdots$, is a basis in the
Hilbert space ${\cal H}_l$. Introducing the notation $\Xi
_n:=\xi_{q^n}(x)$, $n=0,1,2,\cdots$, we find from (3.7) that
$$
J \,\Xi _n = - a q^{-n}(1-bq^{n+1})\,\Xi _{n+1} + q^{-n}(a+1)\,
\Xi _n - q^{-n}(1-q^n)\, \Xi _{n-1} . \eqno (4.3)
$$
As we see, the matrix of the operator $J$ in the basis $\Xi _n$,
$n=0,1,2,\cdots$, is not symmetric, although in the initial basis
$f^l_n$, $n=0,1,2,\cdots$, it was symmetric. The reason is that
the matrix $(a_{mn})$ with entries
$$
a_{mn}:=\beta_m(q^n),\ \ \ \ m,n=0,1,2,\cdots ,
$$
where $\beta_m(q^n)$ are the coefficients (3.4) in the expansion
$\xi _{q^n}(x)=\sum _m \,\beta_m(q^n)f^l_n(x)$, is not unitary. It
is equivalent to the statement that the basis $\Xi
_n:=\xi_{q^n}(x)$, $n=0,1,2,\cdots$, is not normalized. To
normalize it, one has to multiply $\Xi _n$ by corresponding
numbers $c_n$ (which are not known at this moment). Let $\hat\Xi
_n = c_n\Xi _n$, $n=0,1,2,\cdots$, be a normalized basis. Then the
matrix of the operator $J$ is symmetric in this basis. Since  $J$
has in the basis $\{ \hat\Xi _n\}$ the form
$$
J\, \hat\Xi _n = -c_{n+1}^{-1}c_naq^{-n}(1-bq^{n+1})\, \hat\Xi
_{n+1} + q^{-n}(a+1)\, \hat\Xi _n - c_{n-1}^{-1}c_n q^{-n}(1-q^n)
\,\hat\Xi_{n-1} , \eqno (4.4)
$$
then its symmetricity means that
$$
c_{n+1}^{-1}c_naq^{-n}(1-bq^{n+1})=c_{n}^{-1}c_{n+1}
q^{-n-1}(1-q^{n+1})\,,
$$
that is, $c_{n}/c_{n-1} =\sqrt{aq (1-bq^n)/(1-q^n)}$. Therefore,
$$
c_n= c(aq)^{n/2}\frac{(bq;q)_n^{1/2}}{(q;q)_n^{1/2}},
$$
where $c$ is a constant.

The expansions
$$
\hat\xi _{q^n}(x)\equiv \hat\Xi _n(x)=\sum _m
c_n\beta_m(q^n)f^l_m(x) \eqno (4.5)
$$
connect two orthonormal bases in the representation space ${\cal
H}_l$. This means that the matrix $({\hat a}_{mn})$,
$m,n=0,1,2,\cdots$, with entries
$$
{\hat a}_{mn}=c_n\beta _m(q^n)= c\left( (aq)^{n-m}\,
\frac{(bq;q)_n}{(q;q)_n}\,\frac{(abq,aq;q)_m\,(1-abq^{2m+1})}
{(bq,q;q)_m\, (1-abq) }\right) ^{1/2}\, p_m(q^n ;a,b|q) \eqno(4.6)
$$
is unitary, provided  that the constant $c$ is appropriately
chosen. In order to calculate this constant, we use the relation
$\sum_{m=0}^\infty |{\hat a}_{mn}|^2=1$ for $n=0$. Then this sum
is a multiple of the sum in (4.2) and, consequently,
$$
c=\frac{(aq;q)^{1/2}_\infty}{(abq;q)^{1/2}_\infty} .
$$
Thus the $c_n$ in (4.5) and (4.6) is real and equals to
$$
c_n=\left( \frac{(aq;q)_\infty}{(abq;q)_\infty}
\frac{(bq;q)_n(aq)^n}{(q;q)_n} \right) ^{1/2} .
$$

The matrix $({\hat a}_{mn})$ is orthogonal, that is,
$$
\sum _n {\hat a}_{mn}{\hat a}_{m'n}=\delta_{mm'},\ \ \ \ \sum _m
{\hat a}_{mn}{\hat a}_{mn'}=\delta_{nn'} . \eqno (4.7)
$$
Substituting into the first sum over $n$ in (4.7) the expressions
for ${\hat a}_{mn}$, we obtain the identity
$$
\sum_{n=0}^\infty
\frac{(bq;q)_n(aq)^n}{(q;q)_n}\,p_m(q^n;a,b|q)\,p_{m'}(q^n;a,b|q)
$$  $$
=\frac{(abq^2;q)_\infty}{(aq;q)_\infty}
\frac{(1-abq)(aq)^m\,(bq,q;q)_m} {(1-abq^{2m+1})\,(abq,aq;q)_m}\,
\delta_{mm'}\,, \eqno (4.8)
$$
which must yield the orthogonality relation for the little
$q$-Jacobi polynomials. An only gap, which appears here, is the
following. We have assumed that the points $q^n$,
$n=0,1,2,\cdots$, exhaust the whole spectrum of the operator
$I_1$. Let us show that this is the case.

Recall that the selfadjoint operator $I_1$ is represented by a
Jacobi matrix in the basis $f^l_n(x)$, $n=0,1,2,\cdots$. According
to the theory of operators of such type (see, for example, [18],
Chapter VII), eigenfunctions $\xi_\lambda$ of $I_1$ are expanded
into series in the monomials $f^l_n(x)$, $n=0,1,2,\cdots$, with
coefficients, which are polynomials in $\lambda$. These
polynomials are orthogonal with respect to some positive measure
$d\mu (\lambda)$ (moreover, for selfadjoint operators this measure
is unique). The set (a subset of ${\Bbb R}$), on which the
polynomials are orthogonal, coincides with the spectrum of the
operator under consideration and the spectrum is simple. Let us
apply these assertions to the operator $I_1$.

We have found that the spectrum of $I_1$ contains the points
$q^n$, $n=0,1,2,\cdots$. If the operator $I_1$ would have other
spectral points $x$, then on the left-hand side of (4.8) there
would be other summands $\mu_{x_k}\,
p_m({x_k};a,b|q)\,p_{m'}({x_k};a,b|q)$, corresponding to these
additional points. Let us show that these additional summands do
not appear. To this end we set $m=m'=0$ in the relation (4.8) with
the additional summands. Since $p_0(x;a,b|q)=1$, we have the
equality
$$
\sum_{n=0}^\infty \frac{(bq;q)_n(aq)^n}{(q;q)_n} + \sum_k
\mu_{x_k} =\frac{(abq^2;q)_\infty}{(aq;q)_\infty} .
$$
According to the $q$-binomial theorem (see formula (1.3.2) in
[15]), we have
$$
\sum_{n=0}^\infty\,\frac{(bq;q)_n\,(aq)^n}{(q;q)_n} =
\frac{(abq^2;q)_\infty}{(aq;q)_\infty}. \eqno (4.9)
$$
Hence, $\sum_k \mu_{x_k} =0$ and all $\mu_{x_k}$ disappear. This
means that additional summands do not appear in (4.8) and it does
represent the orthogonality relation for the little $q$-Jacobi
polynomials.

By using the operators $I_1$ and $J$, which form a Leonard pair,
we thus derived the orthogonality relation for little $q$-Jacobi
polynomials.

The orthogonality relation for the little $q$-Jacobi polynomials
is given by formula (4.8). Due to this orthogonality, we arrive at
the following statement: {\it The spectrum of the operator $I_1$
coincides with the set of points $q^{n}$, $n=0,1,2,\cdots$. The
spectrum is simple and has one accumulation point at 0.}
\bigskip

\noindent {\bf 5. Dual little $q$-Jacobi polynomials}
\bigskip

Now we consider the second identity in (4.7), which gives the
orthogonality relation for the matrix elements ${\hat a}_{mn}$,
considered as functions of $m$. Up to multiplicative factors these
functions coincide with the functions
$$
F_n(x;a,b|q)={}_2\phi_1 (x,abq/x;\; aq;\; q,q^{n+1}), \eqno (5.1)
$$
considered on the set $x\in \{ q^{-m}\, |\, m=0,1,2,\cdots \}$.
Consequently,
$$
{\hat a}_{mn}= \left( \frac{(aq;q)_\infty}{(abq;q)_\infty}\,
\frac{(bq;q)_n}{(q;q)_n}\,(aq)^{n-m}\,\frac{(abq,aq;q)_m\,(1-abq^{2m+1})}
{(bq,q;q)_m}\right)^{1/2}\, F_n(q^{-m} ;a,b|q)
$$
and the second identity in (4.7) gives the orthogonality relation
for the functions (5.1):
$$
\sum_{m=0}^\infty \,\frac{(1-abq^{2m+1})\,(abq,aq;q)_m}
{(1-abq)(aq)^m\,(bq,q;q)_m}\,
F_n(q^{-m};a,b|q)\,F_{n'}(q^{-m};a,b|q)
$$  $$
= \frac{(abq^2;q)_\infty}{(aq;q)_\infty}
\frac{(q;q)_n(aq)^{-n}}{(bq;q)_n} \, \delta_{nn'}. \eqno (5.2)
$$

The functions $F_n(x;a,b|q)$ can be represented in another form.
Indeed, one can use the relation (III.8) of Appendix III in [15]
in order to obtain that
$$
F_{n}(q^{-m} ;a,b|q)=\frac{(b^{-1}q^{-m};q)_m}{(aq;q)_m}
(abq^{m+1})^m  {}_3\phi_1 (q^{-m},abq^{m+1},q^{-n};\; bq; \;
q,q^n/a)
$$  $$
=\frac{(-1)^m\,(bq;q)_m}{(aq;q)_m}\, a^m \,q^{m(m+1)/2} {}_3\phi_1
(q^{-m},abq^{m+1},q^{-n};\; bq; \; q,q^n/a).  \eqno (5.3)
$$
The basic hypergeometric function ${}_3\phi_1$ in (5.3) is a
polynomial of degree $n$ in the variable $\mu(m): = q^{-m}+
ab\,q^{m+1}$, which represents a $q$-quadratic lattice; we denote
it
$$
d_n(\mu (m); a,b|q):= {}_3\phi_1(q^{-m},ab\,q^{m+1},q^{-n};\;
bq;\; q,q^n/a)\,.  \eqno (5.4)
$$
Then formula (5.2) yields the orthogonality relation
$$
\sum_{m=0}^\infty \,\frac{(1-abq^{2m+1})(abq,bq;q)_m}
{(1-abq)(aq,q;q)_m}\,a^m\, q^{m^2}\, d_n(\mu(m))\, d_{n'}(\mu (m))
$$  $$
= \frac{(abq^2;q)_\infty}{(aq;q)_\infty}
\,\frac{(q;q)_n(aq)^{-n}}{(bq;q)_n}\, \delta_{nn'} \eqno(5.5)
$$
for the polynomials (5.4). We call the polynomials $d_n(\mu (m);
a,b|q)$ {\it dual little $q$-Jacobi polynomials}. Note that these
polynomials can be expressed in terms of the Al-Salam--Chihara
polynomials
$$
Q_n(x;a,b|q)=\frac{(ab;q)_n}{a^n} \left.
{}_3\phi_2\left({q^{-n},az,az^{-1} \atop ab, 0} \right| q,q
\right) , \ \ \ \ x=\frac 12 (z+z^{-1}),
$$
with the parameter $q>1$. An explicit relation between them is
$$
d_n(\mu(x); \beta /\alpha,\, 1/\alpha \beta q \, |\,
q)=\frac{q^{n(n-1)/2}}{(-\beta)^n(1/\alpha\beta ;q)_n}\,
Q_n(\alpha\mu(x)/2; \alpha,\beta | q^{-1}).
$$
Ch. Berg and M. E. H. Ismail studied this type of
Al-Salam--Chihara polynomials in [24] and derived continuous
complex orthogonality measures for them. But [24] does not contain
any discussion of the duality of this family of polynomials with
respect to little $q$-Jacobi polynomials.

A recurrence relation for the polynomials $d_n(\mu (m);a,b|q)$ is
derived from formula (3.6). It has the form
$$
(q^{-m}+ ab q^{m+1})\,d_n(\mu (m))= -\, a\,
q^{-n}(1-bq^{n+1})\,d_{n+1}(\mu (m))
$$   $$
+\, q^{-n}(1+a)\,d_{n}(\mu (m)) - q^{-n}(1-q^{n})\, d_{n-1}(\mu
(m)),
$$
where $d_{n}(\mu (m))\equiv d_{n}(\mu(m);a,b|q)$. Comparing this
relation with the recurrence relation (3.69) in [25], we see that
the polynomials (5.4) are multiple to the polynomials (3.67) in
[25]. Moreover, if one takes into account this multiplicative
factor, the orthogonality relation (5.5) for polynomials (5.4)
turns into relation (3.82) for the polynomials (3.67) in [25],
although the derivation of the orthogonality relation in [25] is
more complicated than our derivation of (5.5). The authors of [25]
do not give an explicit form of their polynomials in the form
similar to (5.4). Concerning the polynomials (3.67) in [25] see
also [26].

Let ${\frak l}^2$ be the Hilbert space of functions on the set
$m=0,1,2,\cdots$ with the scalar product
$$
\langle f_1,f_2\rangle = \sum _{m=0}^\infty\,
\frac{(1-abq^{2m+1})\, (abq,bq;q)_m}{(1-abq)\, (aq,q;q)_m}\,
a^m\,q^{m^2}\,f_1(m)\,\overline{f_2(m)} .
$$
The polynomials (5.4) are in one-to-one correspondence with the
columns of the unitary matrix $({\hat a}_{mn})$ and the
orthogonality relation (5.5) is equivalent to the orthogonality of
these columns. Due to (4.7) the columns of the matrix $({\hat
a}_{mn})$ form an orthonormal basis in the Hilbert space of
sequences ${\bf a}=\{ a_n\, |\, n=0,1,2,\cdots \}$ with the scalar
product $\langle {\bf a},{\bf a}'\rangle=\sum_n a_na'_n$. For this
reason, the set of polynomials $d_n(\mu (m);a,b|q)$,
$n=0,1,2,\cdots$, form an orthogonal basis in the Hilbert space
${\frak l}^2$. This means that {\it the point measure in (5.5) is
extremal for the dual little $q$-Jacobi polynomials} $d_n(\mu
(m);a,b|q)$ (for the definition of an extremal orthogonality
measure see, for example, [27]).
\bigskip

\noindent{\bf 6. Pair of operators $(I_2,J)$}
\bigskip

Let $b$ and $c$ be real numbers such that $0<b<q^{-1}$ and $c<0$.
We consider the operator
$$
I_2 := \alpha\,q^{J_0/4}(J_+\, A +A\, J_-)\,q^{J_0/4} - B \eqno
(6.1)
$$
of the representation $T^+_l$, where $\alpha
=(-ac/q)^{1/2}\,(1-q)$, $a = q^{2l-1}$, and $A$ and $B$ are
operators of the form
$$
A=\frac{q^{(J_0-l+2)/2} \sqrt{(1-bq^{J_0-l+1})(1-abq^{J_0-l+1})
(1-cq^{J_0-l+1})(1-abc^{-1}q^{J_0-l+1})}}
{(1-abq^{2J_0-2l+2})\sqrt{(1-abq^{2J_0-2l+1})(1-abq^{2J_0-2l+3})}}\,,
$$   $$
B= \frac{(1-aq^{J_0-l+1})(1-abq^{J_0-l+1})(1-cq^{J_0-l+1})}
{(1-abq^{2J_0-2l+1)} (1-abq^{2J_0-2l+2})}
$$  $$
-
ac\,q^{J_0-l+1}\,\frac{(1-q^{J_0-l})(1-bq^{J_0-l})(1-abc^{-1}q^{J_0-l})}
{(1-abq^{2J_0-2l})(1-abq^{2J_0-2l+1})} - 1.
$$

We have the following formula for the symmetric operator $I_2$:
$$
I_2\, f^l_n=a_nf^l_{n+1}+a_{n-1}f^l_{n-1}-b_nf^l_n , \eqno (6.2)
$$
where
$$
a_{n-1}=(-acq^{n+1})^{1/2}\,\frac{ \sqrt{(1-q^{n})(1-aq^{n})
(1-bq^{n})(1-abq^{n})\,(1-cq^{n})(1-abc^{-1}q^{n})}}
{(1-abq^{2n})\,\sqrt{(1-abq^{2n-1})(1-abq^{2n+1})}} ,
$$  $$
b_n= \frac{(1-aq^{n+1})(1-abq^{n+1})(1-cq^{n+1})}
{(1-abq^{2n+1})(1-abq^{2n+2})}-acq^{n+1}
\frac{(1-q^{n})(1-bq^{n})(1-abq^{n}/c)}{(1-abq^{2n})
(1-abq^{2n+1})} - 1\,,
$$
with $a=q^{2l-1}$. Recall that $l$ takes any positive value, so
$a$ can be any positive number such that $a<q^{-1}$. Since the
$q^{\pm J_0}$ are symmetric operators, the operator $I_2$ is also
symmetric.

The operator $I_2$ is bounded. Therefore, we assume that it is
defined on the whole representation space ${\cal H}_l$. This means
that $I_2$ is a selfadjoint operator. Actually, $I_2$ is a trace
class operator. To show this we note that for the coefficients
$a_n$ and $b_n$ from (6.2) one obtains that
$$
a_{n+1}/a_n \to q^{1/2},\ \ b_{n+1}/b_n \to q \ \ \ {\rm when}\ \
\ n\to \infty .
$$
Therefore, $\sum _n (2a_n+b_n)< \infty$ and this means that $I_2$
is a trace class operator. Thus, the spectrum of $I_2$ is simple
(since it is representable by a Jacobi matrix with $a_n\ne 0$),
discrete and have a single accumulation point at 0.

To find eigenfunctions $\psi_\lambda (x)$ of the operator $I_2$,
$I_2 \psi_\lambda (x)=\lambda \psi_\lambda (x)$, we set
$$
\psi_\lambda (x)=\sum _{n=0}^\infty \beta_n(\lambda)f^l_n (x).
$$
Acting by the operator $I_2$ on both sides of this relation, one
derives that
$$
\sum_n \beta_n(\lambda) (a_nf^l_{n+1}+a_{n-1}f^l_{n-1}-b_nf^l_n)
=\lambda \sum \beta_n(\lambda)f^l_n ,
$$
where $a_n$ and $b_n$ are the same as in (6.2). Collecting in this
identity factors, which multiply $f^l_n$ with fixed $n$, we arrive
at the recurrence relation for the coefficients
$\beta_n(\lambda)$:
$$
a_n\beta_{n+1}(\lambda) +a_{n-1}\beta_{n-1}(\lambda)-
b_n\beta_{n}(\lambda) = \lambda \beta_{n}(\lambda).
$$
Making the substitution
$$
\beta_{n}(\lambda)=\left(\frac{(abq,aq,cq;q)_n\,
(1-abq^{2n+1})}{(abq/c,bq,q;q)_n\,
(1-abq)(-ac)^n}\right)^{1/2}q^{-n(n+3)/4} \,\beta'_{n}(\lambda)
$$
we reduce this relation to the following one
$$
A_n \beta'_{n+1}(\lambda)+ C_n \beta'_{n-1}(\lambda)
-(A_n+C_n-1)p'_{n}(\lambda)=\lambda \beta'_{n}(\lambda)
$$
with
$$
A_n=\frac{(1-aq^{n+1})(1-cq^{n+1})(1-abq^{n+1})}
{(1-abq^{2n+1})\,(1-abq^{2n+2})},
$$  $$
C_n=\frac{-acq^{n+1}(1-q^{n})(1-bq^{n})(1-abc^{-1}q^{n})}
{(1-abq^{2n})\,(1-abq^{2n+1})}.
$$
It is the recurrence relation for the big $q$-Jacobi polynomials
$$
P_n(\lambda ;a,b,c;q):= {}_3\phi_2 (q^{-n}, abq^{n+1}, \lambda ;\;
aq,cq; \; q,q )  \eqno (6.3)
$$
introduced by G. E. Andrews and R. Askey [28] (see also formula
(7.3.10) in [15]). Therefore, $\beta'_n(\lambda )=P_n(\lambda
;a,b,c;q)$ and
$$
\beta_n(\lambda )=\left(\frac{(abq,aq,cq;q)_n\,(1-abq^{2n+1})}
{(abq/c,bq,q;q)_n\, (1-abq) (-ac)^n} \right)^{1/2}q^{-n(n+3)/4}\,
P_n(\lambda ;a,b,c;q). \eqno (6.4)
$$
For the eigenfunctions $\psi _\lambda(x)$ we have the expansion
$$
\psi _\lambda(x)=\sum_{n=0}^\infty \left(\frac{(abq,aq,cq;q)_n\,
(1-abq^{2n+1})}{(abq/c,bq,q;q)_n\,
(1-abq)(-ac)^n}\right)^{1/2}q^{-n(n+3)/4}\, P_n(\lambda
;a,b,c;q)\, f^l_n(x)
$$   $$
=\sum_{n=0}^\infty \,a^{-n/4}\, \frac{(aq;q)_n}{(q;q)_n}
\left(\frac{(abq,cq;q)_n\,(1-abq^{2n+1})} {(abq/c,bq;q)_n\,
(1-abq) (-ac)^n} \right)^{1/2}q^{-n(n+3)/4}\, P_n(\lambda
;a,b,c;q)\, x^n. \eqno (6.5)
$$
Since the spectrum of the operator $I_2$ is discrete, only a
discrete set of these functions belongs to the Hilbert space
${\cal H}_l$ and this discrete set determines the spectrum of
$I_2$.

In what follows we intend to study a spectrum of the operator
$I_2$ and to find polynomials, dual to big $q$-Jacobi polynomials.
It can be done with the aid of the operator
$$
J:= q^{-J_0+l} + ab\,q^{J_0-l+1}\equiv \mu (J_0-l),
$$
which has been already used in the previous case in section 3. In
order to determine how this operator acts upon the eigenfunctions
$\psi _\lambda(x)$, one can use the $q$-difference equation
$$
(q^{-n}+abq^{n+1})\, P_n(\lambda)= aq \lambda^{-2}( \lambda
-1)(b\lambda -c)\,P_{n}(q\lambda )
$$
$$
- [\lambda^{-2}acq(1+q)-\lambda^{-1}q(ab+ac+a+c)]\, P_n(\lambda)+
\lambda^{-2}(\lambda-aq)(\lambda-cq)\,P_n(q^{-1}\lambda) ,
\eqno(6.6)
$$
for the big $q$-Jacobi polynomials $P_n(\lambda)\equiv P_n
(\lambda;a,b,c;q)$ (see, for example, formula (3.5.5) in [19]).
Multiply both sides of (6.6) by $d_n\,f^l_n(x)$, where $d_n$ are
the coefficients of $P_n(\lambda ;a,b,c;q)$ in the expression
(6.4) for the coefficients $\beta_n(\lambda)$, and sum over $n$.
Taking into account formula (6.5) and the fact that
$J\,f^l_n(x)=(q^{-n}+ab\,q^{n+1})\,f^l_n(x)$, one obtains the
relation
$$
J\, \psi _{\lambda}(x)= aq\lambda^{-2}( \lambda -1)(b\lambda -c)
\,\psi_{q\lambda}(x)
$$
$$
- [\lambda^{-2}acq(1+q)-\lambda^{-1}q(ab+ac+a+c)]\,\psi
_{\lambda}(x) + \lambda^{-2}(\lambda-aq)(\lambda-cq)\,\psi
_{q^{-1}\lambda}(x) .\eqno (6.7)
$$
It will be shown in the next section that the spectrum of the
operator $I_2$ consists of the points  $aq^n$, $cq^n$,
$n=0,1,2,\cdots$. The matrix $J$ consists of two Jacobi matrices
(one corresponds to the basis elements $aq^n$, $n=0,1,2,\cdots$,
and another to the basis elements $cq^n$, $n=0,1,2,\cdots$). In
this case, the operators $I_2$ and $J$ form some generalization of
a Leonard pair.
\bigskip

\noindent{\bf 7. Spectrum of $I_2$ and orthogonality of big
$q$-Jacobi polynomials}
\bigskip

As in section 4 one can show that for some value of $\lambda$
(which must belong to the spectrum) the last term on the right
side of (6.7) has to vanish. There are two such values of
$\lambda$: $\lambda = aq$ and $\lambda = cq$. Let us show that
both of these points are spectral points of the operator $I_2$.
Observe that, according to (6.3),
$$
P_n(aq ;a,b,c;q):={}_2\phi_1 (q^{-n}, abq^{n+1} ;\; cq; \; q,q )
=\frac{(c/abq^n;q)_n}{(cq;q)_n} (ab)^n\, q^{n(n+1)}.
$$
Therefore, since
$$
(c/abq^n;q)_n=(abq/c;q)_n(-c/ab)^n q^{-n(n+1)/2}\,,
$$
one obtains that
$$
P_n(aq ;a,b,c;q):=\frac{(abq/c;q)_n}{(cq;q)_n}(-c)^nq^{n(n+1)/2}.
\eqno (7.1)
$$
Likewise,
$$
P_n(cq ;a,b,c;q):=\frac{(bq;q)_n}{(aq;q)_n}(-a)^nq^{n(n+1)/2}.
$$
Hence, for the scalar product $\langle \psi_{aq}(x),
\psi_{aq}(x)\rangle$ we have the expression
$$
\sum_{n=0}^\infty\,\frac{(1-abq^{2n+1})\,(abq,aq,cq;q)_n}
{(1-abq)(abq/c,bq,q;q)_n\, (-ac)^n}\,q^{-n(n+3)/2}\,
P^2_n(aq;a,b,c;q)
$$  $$
= \sum_{n=0}^\infty\,\frac{(1-abq^{2n+1})\,(abq/c,abq,aq;q)_n}
{(1-abq)(bq,cq,q;q)_n\,(-a/c)^n}\,q^{n(n-1)/2} =
\frac{(abq^2,c/a;q)_\infty}{(bq,cq;q)_\infty}, \eqno (7.2)
$$
where the relation (A.6) from Appendix has been used. Similarly,
for $\langle \psi_{cq}(x),\psi_{cq}(x) \rangle$ one has the
expression
$$
\sum_{n=0}^\infty\,\frac{(1-abq^{2n+1})(abq,aq,cq;q)_n}
{(1-abq)(-ac)^n\,(abq/c,bq,q;q)_n}\,q^{-n(n+3)/2}\,P^2_n(cq;a,b,c;q)
= \frac{(abq^2,a/c;q)_\infty}{(aq,abq/c;q)_\infty} , \eqno(7.3)
$$
where formula (A.7) from Appendix has been used. Thus, the values
$\lambda=aq$ and $\lambda=cq$ are the spectral points of the
operator $I_2$.

Let us find other spectral points of the operator $I_2$. Setting
$\lambda = aq$ in (6.7), we see that the operator $J$ transforms
$\psi _{aq}(x)$ into a linear combination of the functions
$\psi_{aq^2}(x)$ and $\psi _{aq}(x)$. We have to show that
$\psi_{aq^2}(x)$ belongs to the Hilbert space ${\cal H}_l$, that
is, that
$$
\langle \psi _{aq^2} ,\psi _{aq^2} \rangle = \sum_{n=0}^\infty\,
\frac{(abq,aq,cq;q)_n\,(1-abq^{2n+1})}{(abq/c,bq,q;q)_n\,
(1-abq)(-ac)^n} \,q^{-n(n+3)/2}\, P^2_n(aq^2;a,b,c;q)<\infty .
$$
It is harder to prove this inequality directly, than in the case
of the little $q$-Jacobi polynomials. Therefore we do not want to
embark on the discussion of this point here (for this inequality
actually follows from the orthogonality relation for the big
$q$-Jacobi polynomials). The above inequality shows that $\psi
_{aq^2}(x)$ is an eigenfunction of $I_2$ and the point $aq^2$
belongs to the spectrum of the operator $I_2$. Setting $\lambda =
aq^2$ in (6.7) and acting similarly, one obtains that $\psi
_{aq^3}(x)$ is an eigenfunction of $I_2$ and the point $aq^3$
belongs to the spectrum of $I_2$. Repeating this procedure, one
sees that $\psi _{aq^n}(x)$, $n=1,2,\cdots$, are eigenfunctions of
$I_2$ and the set $aq^n$, $n=1,2,\cdots$, belongs to the spectrum
of $I_2$. Likewise, one concludes that $\psi _{cq^n}(x)$,
$n=1,2,\cdots$, are eigenfunctions of $I_2$ and the set $cq^n$,
$n=1,2,\cdots$, belongs to the spectrum of $I_2$. Note that so far
we do not know whether the operator $I_2$ has other spectral
points or not. In order to solve this problem we shall proceed as
in section 4.

The functions $\psi _{aq^n}(x)$ and $\psi _{cq^n}(x)$,
$n=1,2,\cdots$, are linearly independent elements of the
representation space ${\cal H}_l$. Suppose that $aq^n$ and $cq^n$,
$n=1,2,\cdots$, constitute the whole spectrum of the operator
$I_2$. Then the set of functions $\psi _{aq^n}(x)$ and
$\psi_{cq^n}(x)$, $n=1,2,\cdots$, is a basis of the Hilbert space
${\cal H}_l$. Introducing the notations $\Xi
_n:=\xi_{aq^{n+1}}(x)$ and $\Xi'_n:=\xi _{cq^{n+1}}(x)$,
$n=0,1,2,\cdots$, we find from (6.7) that
$$
J\, \Xi _n = a^{-1}cq^{-2n-1}(1-aq^{n+1})(1-baq^{n+1}/c)\, \Xi
_{n+1} + d_n\,\Xi _n + a^{-1}cq^{-2n}(1-q^n)(1-aq^n/c)\, \Xi
_{n-1} ,
$$ $$
J \,\Xi'_n = c^{-1}aq^{-2n-1}(1-cq^{n+1})(1-bq^{n+1})\, \Xi _{n+1}
+ d'_n\,\Xi _n + c^{-1}aq^{-2n}(1-q^n)(1-cq^n/a)\, \Xi _{n-1} ,
$$
where
$$
d_n= \frac 1a [q^{-2n-1}c(1+q)-q^{-n}(ab+ac+a+c)] ,
$$  $$
d'_n= \frac 1c [q^{-2n-1}a(1+q)-q^{-n}(ab+ac+a+c)] .
$$

As we see, the matrix of the operator $J$ in the basis $\Xi _n
=\xi_{aq^{n+1}}(x)$, $\Xi'_n=\xi _{cq^{n+1}}(x)$,
$n=0,1,2,\cdots$, is not symmetric, although in the initial basis
$f^l_n$, $n=0,1,2,\cdots$, it was symmetric. The reason is that
the matrix $M:=(a_{mn}\ \;  a'_{mn})$ with entries
$$
a_{mn}:=\beta_m(aq^n),\ \ \  a'_{mn}:=\beta_m(cq^n),\ \ \
m,n=0,1,2,\cdots ,
$$
where $\beta_m(dq^n)$, $d=a,c$, are coefficients (6.4) in the
expansion $\psi _{dq^n}(x)=\sum _m \,\beta_m(dq^n)\,f^l_n(x)$ (see
above), is not unitary. This matrix $M$ is formed by adding the
columns of the matrix $(a'_{mn})$ to the columns of the matrix
$(a_{mn})$ from the right, that is,
$$
M =\left( \matrix{
 a_{11}&\cdots &a_{1k} &\cdots &a'_{11} &\cdots &a'_{1l}&\cdots\cr
 a_{21}&\cdots &a_{2k} &\cdots &a'_{21} &\cdots &a'_{2l}&\cdots\cr
 \cdots&\cdots &\cdots &\cdots &\cdots &\cdots  &\cdots &\cdots\cr
 a_{j1}&\cdots &a_{jk} &\cdots &a'_{j1} &\cdots &a'_{jl}&\cdots\cr
 \cdots&\cdots &\cdots &\cdots &\cdots &\cdots  &\cdots
 &\cdots\cr}
\right) .
$$
It maps the basis $\{ f^l_n\}$ into the basis $\{\psi_{aq^{n+1}},
\psi _{cq^{n+1}} \}$ in the representation space. The nonunitarity
of the matrix $M$ is equivalent to the statement that the basis
$\Xi _n:=\xi _{aq^{n+1}}(x)$, $\Xi _n:=\xi _{cq^{n+1}}(x)$,
$n=0,1,2,\cdots$, is not normalized. In order to normalize it we
have to multiply $\Xi _n$ by appropriate numbers $c_n$ and
$\Xi'_n$ by numbers $c'_n$. Let $\hat\Xi _n = c_n\Xi _n$,
$\hat\Xi'_n =c'_n\Xi_n$, $n=0,1,2,\cdots$, be a normalized basis.
Then the operator $J$ is symmetric in this basis and has the form
$$
J\,\hat\Xi _n = c_{n+1}^{-1}c_na^{-1}cq^{-2n-1}(1-aq^{n+1})
(1-abq^{n+1}/c)\,\hat\Xi_{n+1} + d_n\, \hat\Xi _n
$$  $$
+ c_{n-1}^{-1}c_n a^{-1}cq^{-2n}(1-aq^{n}/c)(1-q^n)\,\hat\Xi
_{n-1}, \eqno (7.4)
$$
$$
J\,\hat\Xi'_n = {c'}_{n+1}^{-1}{c'}_nc^{-1}aq^{-2n-1}(1-bq^{n+1})
(1-cq^{n+1})\,\hat\Xi _{n+1} + d'_n \,\hat\Xi _n
$$  $$
-{c'}_{n-1}^{-1}{c'}_n c^{-1}aq^{-2n}
(1-cq^{n}/a)(1-q^n)\,\hat\Xi_{n-1} , \eqno (7.5)
$$
The symmetricity of the matrix of the operator $J$ in the basis
$\{ \hat\Xi _n,\hat\Xi'_n \}$ means that
$$
c_{n+1}^{-1}c_n q^{-2n-1}(1-aq^{n+1})(1-abq^{n+1}/c)
=c_{n}^{-1}c_{n+1} q^{-2n-2} (1-aq^{n+1}/c)(1-q^{n+1}),
$$  $$
{c'}_{n+1}^{-1}{c'}_n q^{-2n-1}(1-bq^{n+1})(1-cq^{n+1})
={c'}_{n}^{-1}{c'}_{n+1} q^{-2n-2} (1-cq^{n+1}/a)(1-q^{n+1}).
$$
that is,
$$
\frac{c_{n}}{c_{n-1}} =\sqrt{q\frac{(1-aq^n)(1-abq^n/c)}
{(1-q^n)(1-aq^n/c)}} , \ \ \ \frac{c'_{n}}{c'_{n-1}}
=\sqrt{q\frac{(1-cq^n)(1-bq^n)} {(1-q^n)(1-cq^n/a)}} .
$$
Thus,
$$
c_n= C\left(
q^{n}\frac{(abq/c,aq;q)_n}{(aq/c,q;q)_n}\right)^{1/2},\ \ \ c'_n=
C'\left(q^{n}\frac{(bq,cq;q)_n}{(cq/a,q;q)_n}\right)^{1/2},
$$
where $C$ and $C'$ are some constants.

Therefore, in the expansions
$$
\hat\psi _{aq^n}(x)\equiv \hat\Xi _n(x)= \sum _m
\,c_n\,\beta_m(aq^n)\,f^l_m(x)=\sum _m \,{\hat a}_{mn}\, f^l_m(x),
\eqno (7.6)
$$
$$
\hat\psi _{cq^n}(x)\equiv \hat\Xi _n(x) = \sum _m
\,c'_n\,\beta_m(cq^n)\,f^l_m(x)=\sum _m \,{\hat a}'_{mn}\,
f^l_m(x), \eqno (7.7)
$$
the matrix ${\hat M}:=({\hat a}_{mn}\ {\hat a}'_{mn})$ with
entries
$$
{\hat a}_{mn}= c_n\,\beta _m(aq^n) = C
\left(q^{n}\frac{(abq/c,aq;q)_n}{(aq/c,q;q)_n}
\frac{(abq,aq,cq;q)_m\,(1-abq^{2m+1})}{(abq/c,bq,q;q)_m\,
(1-abq)(-ac)^m }\right)^{1/2}
$$
$$
\times \, q^{-m(m+3)/4}\,P_m(aq^{n+1} ;a,b,c;q), \eqno(7.8)
$$  $$
{\hat a}'_{mn}= c_n\,\beta _m(cq^n)= C'
\left(q^{n}\frac{(bq,cq;q)_n}{((cq/a,q;q)_n}\,
\frac{(abq,aq,cq;q)_m\,(1-abq^{2m+1})}{(abq/c,bq,q;q)_m
\,(1-abq)(-ac)^m}\right)^{1/2} $$
$$ \times \,q^{-m(m+3)/4}\,P_m(cq^{n+1} ;a,b,c;q)\,,
\eqno(7.9)
$$
is unitary, provided that the constants $C$ and $C'$ are
appropriately chosen. In order to calculate these constants, one
can use the relations
$$
\sum _{m=0}^\infty |{\hat a}_{mn}|^2=1\,,\ \ \ \ \sum
_{m=0}^\infty |{\hat a}'_{mn}|^2=1 \,,
$$
for $n=0$. Then these sums are multiples of the sums in (7.2) and
(7.3), so we find that
$$
C=\frac{(bq,cq;q)^{1/2}_\infty}{(abq^2, c/a;q)^{1/2}_\infty} ,\ \
\ C'=\frac{(aq,abq/c;q)^{1/2}_\infty}{(abq^2,
a/c;q)^{1/2}_\infty}.\eqno (7.10)
$$
The coefficients $c_n$ and $c'_n$ in (7.6)--(7.9) are thus real
and equal to
$$
c_n= \left( \frac{(bq,cq;q)_\infty}{(abq^2, c/a;q)_\infty}
\frac{(abq/c,aq;q)_n\, q^n}{(aq/c,q;q)_n} \right)^{1/2}, \ \ \ \
c'_n= \left( \frac{(aq,abq/c;q)_\infty}{(abq^2, a/c;q)_\infty}
\frac{(bq,cq;q)_n\, q^n}{(cq/a,q;q)_n} \right)^{1/2}.
$$

The orthogonality of the matrix ${\hat M}\equiv ({\hat a}_{mn}\
{\hat a}'_{mn})$ means that
$$
\sum _m {\hat a}_{mn}{\hat a}_{mn'}=\delta_{nn'},\ \ \ \sum _m
{\hat a}'_{mn}{\hat a}'_{mn'}=\delta_{nn'},\ \ \ \sum _m {\hat
a}_{mn}{\hat a}'_{mn'}=0,       \eqno (7.11)
$$  $$
\sum _n ({\hat a}_{mn}{\hat a}_{m'n}+ {\hat a}'_{mn} {\hat
a}'_{m'n} ) =\delta_{mm'} .               \eqno (7.12)
$$
Substituting the expressions for ${\hat a}_{mn}$ and ${\hat
a}'_{mn}$ into (7.12), one obtains the relation
$$
\frac{(bq,
cq;q)_\infty}{(abq^2,c/a;q)_\infty}\,\sum_{n=0}^\infty\,
\frac{(aq,abq/c;q)_n
q^n}{(aq/c,q;q)_n}P_m(aq^{n+1})P_{m'}(aq^{n+1})
$$  $$
+\frac{(aq,
abq/c;q)_\infty}{(abq^2,a/c;q)_\infty}\,\sum_{n=0}^\infty\,
\frac{(bq,cq;q)_n q^n}{(cq/a,q;q)_n}P_m(cq^{n+1})P_{m'}(cq^{n+1})
$$  $$
=\frac{(1-abq)(bq,abq/c,q;q)_m}{(1-abq^{2m+1})(aq,abq,cq;q)_m}
 (-ac)^m q^{m(m+3)/2}\,\delta_{mm'}\, .        \eqno (7.13)
$$
This identity must give an orthogonality relation for the big
$q$-Jacobi polynomials $P_m(y)\equiv P_m(y;a,b,c;q)$. An only gap,
which appears here, is the following. We have assumed that the
points $aq^n$ and $cq^n$, $n=0,1,2,\cdots$, exhaust the whole
spectrum of the operator $I_2$. As in the case of the operator
$I_1$ in section 4, if the operator $I_2$ would have other
spectral points $x_k$, then on the left-hand side of (7.13) would
appear other summands $\mu_{x_k}
P_m({x_k};a,b,c;q)P_{m'}({x_k};a,b,c;q)$, which correspond to
these additional points. Let us show that these additional
summands do not appear. To this end we set $m=m'=0$ in the
relation (7.13) with the additional summands. This results in the
equality
$$
\frac{(bq,
cq;q)_\infty}{(abq^2,c/a;q)_\infty}\,\sum_{n=0}^\infty\,
\frac{(aq,abq/c;q)_n q^n}{(aq/c,q;q)_n}
$$  $$
+\frac{(aq, abq/c;q)_\infty}{(abq^2,a/c;q)_\infty}\,
\sum_{n=0}^\infty \,\frac{(bq,cq;q)_n q^n}{(cq/a,q;q)_n} +
\sum_k\,\mu_{x_k} =1.                           \eqno (7.14)
$$
In order to show that $\sum_k \mu_{x_k} = 0$, take into account
the relation
$$
\frac{(Aq/C, Bq/C;q)_\infty}{(q/C,ABq/C;q)_\infty}\, {}_2\phi_1
(A,B;C;\, q,q)
$$   $$
\frac{(A, B;q)_\infty}{(C/q,ABq/C;q)_\infty}\, {}_2\phi_1
(Aq/C,Bq/C;q^2/C;\, q,q) = 1
$$
(see formula (2.10.13) in [15]). Putting here $A=aq$, $B=abq/c$
and $C=aq/c$, we obtain relation (7.14) without the summand
$\sum_k \mu_{x_k}$. Therefore, in (7.14) the sum $\sum_k
\mu_{x_k}$ does really vanish and formula (7.13) gives an
orthogonality relation for big $q$-Jacobi polynomials.

By using the operators $I_2$ and $J$, we thus derived the
orthogonality relation for big $q$-Jacobi polynomials.

The orthogonality relation (7.13) for big $q$-Jacobi polynomials
enables us to formulate the following statement: {\it The spectrum
of the operator $I_2$ coincides with the set of points $aq^{n+1}$
and $cq^{n+1}$, $n=0,1,2,\cdots$. The spectrum is simple and has
one accumulation point at 0.}
\bigskip

\noindent {\bf 8. Dual big $q$-Jacobi polynomials}
\bigskip

Now we consider the relations (7.11). They give the orthogonality
relation for the set of matrix elements ${\hat a}_{mn}$ and ${\hat
a}'_{mn}$, viewed as functions of $m$. Up to multiplicative
factors, they coincide with the functions
 $$
F_n(x;a,b,c;q):={}_3\phi_2 (x,abq/x, aq^{n+1};\; aq, cq;\; q,q),\
\ n=0,1,2,\cdots,  \eqno (8.1)
$$   $$
F'_n(x;a,b,c;q):={}_3\phi_2 (x,abq/x, cq^{n+1};\; aq, cq;\;
q,q)\equiv F_n(x;c,ab/c,a),\ \ n=0,1,2,\cdots ,  \eqno (8.2)
$$
considered on the corresponding sets of points. Namely, we have
$$
{\hat a}_{mn}\equiv {\hat
a}_{mn}(a,b,c)=C\,\left(q^{n}\frac{(abq/c,aq;q)_n}{(aq/c,q;q)_n}
\frac{(abq,aq,cq;q)_m\,(1-abq^{2m+1})}{(abq/c,bq,q;q)_m\,
(1-abq)(-ac)^m}\right)^{1/2}
$$  $$
\times q^{-m(m+3)/4}\,F_n(q^{-m};a,b,c;q),   \eqno (8.3)
$$   $$
{\hat a}'_{mn}\equiv {\hat
a}'_{mn}(a,b,c)=C'\,\left(q^{n}\frac{(bq,cq;q)_n}{(cq/a,q;q)_n}
\,\frac{(abq,aq,cq;q)_m\,(1-abq^{2m+1})}{(abq/c,bq,q;q)_m\,
(1-abq)(-ac)^m}\right)^{1/2}
$$  $$
\times q^{-m(m+3)/4}\,F'_n(q^{-m} ;a,b,c;q)\equiv {\hat
a}_{mn}(c,ab/c,a), \eqno (8.4)
$$
where $C$ and $C'$ are given by formulas (7.10). The relations
(7.11) lead to the following orthogonality relation for the
functions (8.1) and (8.2):
$$
\frac{(bq, cq;q)_\infty}{(abq^2,c/a;q)_\infty}\,\sum_{m=0}^\infty
\, \rho (m) F_n(q^{-m};a,b,c;q)\,F_{n'}(q^{-m};a,b,c;q)
=\frac{(aq/c,q;q)_n}{(aq,abq/c;q)_n q^n}\delta_{nn'}, \eqno (8.5)
$$   $$
\frac{(aq, abq/c;q)_\infty}{(abq^2,a/c;q)_\infty}\,
\sum_{m=0}^\infty\, \rho (m)
F'_n(q^{-m};a,b,c;q)\,F'_{n'}(q^{-m};a,b,c;q)
=\frac{(cq/a,q;q)_n}{(bq,cq;q)_n q^n}\,\delta_{nn'},  \eqno (8.6)
$$   $$
\sum_{m=0}^\infty\, \rho (m)
F_n(q^{-m};a,b,c;q)\,F'_{n'}(q^{-m};a,b,c;q)=0,  \eqno (8.7)
$$
where
$$
\rho (m):=\frac{(1-abq^{2m+1})(aq,abq,cq;q)_m}{(1-abq)
(bq,abq/c,q;q)_m\,(-ac)^m}\, q^{-m(m+3)/2}.
$$

There is another form for the functions $F_n(q^{-m};a,b,c;q)$ and
$F'_n(q^{-m};a,b,c;q)$. Indeed, one can use the relation (III.12)
of Appendix III in [15] to obtain that
$$
F_n(q^{-m};a,b,c;q)=\frac{(cq^{-m}/ab;q)_m}{(cq;q)_m}
(abq^{m+1})^m \, {}_3\phi_2\left( \left. {q^{-m},abq^{m+1},q^{-n}
   \atop aq, abq/c} \right| q,aq^{n+1}/c  \right)
$$  $$
=\frac{(abq/c;q)_m}{(cq;q)_m} (-c)^{m}q^{m(m+1)/2} \,
{}_3\phi_2\left( \left. {q^{-m},abq^{m+1},q^{-n}
   \atop aq, abq/c} \right| q,aq^{n+1}/c  \right)
$$
and
$$
F'_n(q^{-m};a,b,c;q) =\frac{(bq;q)_m}{(aq;q)_m}
(-a)^{m}q^{m(m+1)/2}\, {}_3\phi_2\left( \left.
{q^{-m},abq^{m+1},q^{-n}
   \atop bq, cq} \right| q,cq^{n+1}/a  \right) .
$$
The basic hypergeometric functions ${}_3\phi_2$ in these formulas
are polynomials in $\mu (m):=q^{-m}+ab\,q^{m+1}$. So if we
introduce the notation
$$
D_n(\mu (m); a,b,c|q):=\left.
{}_3\phi_2\left({q^{-m},abq^{m+1},q^{-n}
  \atop aq, abq/c} \right| q,aq^{n+1}/c \right) ,       \eqno(8.8)
$$
then
$$
F_n(q^{-m};a,b,c;q)=\frac{(abq/c;q)_m}{(cq;q)_m}
(-c)^{m}q^{m(m+1)/2} D_n(\mu (m); a,b,c|q),
$$   $$
F'_n(q^{-m};a,b,c;q)=\frac{(bq;q)_m}{(aq;q)_m}
(-a)^{m}q^{m(m+1)/2} D_n(\mu (m); b,a,ab/c|q).
$$

Formula (8.5) directly leads to the orthogonality relation for the
polynomials $D_n(\mu(m))\equiv D_n(\mu (m); a,b,c|q)$:
$$
\sum_{m=0}^\infty
\frac{(1-abq^{2m+1})(aq,abq,abq/c;q)_m}{(1-abq)(bq,cq,q;q)_m}
\,(-c/a)^m\, q^{m(m-1)/2}\,D_n(\mu (m))\,D_{n'}(\mu (m))
$$   $$
=\frac{(abq^2,c/a;q)_\infty}{(bq,cq;q)_\infty} \frac{(aq/c,q;q)_n}
{(aq,abq/c;q)_nq^n} \delta_{nn'}.  \eqno (8.9)
$$
From (8.6) one obtains the orthogonality relation for the
polynomials $D_n(\mu (m);b,a,ab/c|q)$ (which follows also from the
relation (8.9) by interchanging $a$ and $b$ and replacing $c$ by
$ab/c$).

We call the polynomials $D_n(\mu (m);a,b,c|q)$ {\it dual big
$q$-Jacobi polynomials}. It is natural to ask whether they can be
identified with some known and thoroughly studied set of
polynomials. The answer is: they can be obtained from the
$q$-Racah polynomials $R_n(\mu (x);a,b,c,d|q)$ of Askey and Wilson
[29] by setting $a=q^{-N-1}$ and sending $N\to \infty$, that is,
$$
D_n(\mu (x);a,b,c|q)=\lim_{N\to \infty} R_n( \mu (x); q^{-N-1},
a/c,a,b|q).   \eqno (8.10)
$$
Observe that the orthogonality relation (8.9) can be also derived
from formula (4.16) in [30]. But the derivation of this formula
(4.16) is rather complicated.

It is worth noting here that in the limit as $c\to 0$ the dual big
$q$-Jacobi polynomials $D_n(\mu (x);a,b,c|q)$ coincide with the
dual little $q$-Jacobi polynomials $d_n(\mu (x);b,a|q)$, defined
in section 5. The dual little $q$-Jacobi polynomials $d_n(\mu
(x);a,b|q)$ reduce, in turn, to the Al-Salam--Carlitz II
polynomials $V_n^{(a)}(s;q)$ on the $q$-linear lattice $s=q^{-x}$
(see [19], p. 114) in the case when the parameter $b$ vanishes,
that is,
$$
d_n(\mu (x);a,0|q)=´ç_2\phi_0(q^{-n},\, q^{-x}; - ;\; q,q^n/a)=
(-a)^{-n}q^{n(n-1)/2}\, V_n^{(a)}(q^{-x};q). \eqno (8.11)
$$
This means that we have now a complete chain of reductions
$$
R_n(\mu (x);a,b,c,d|q)\; {\mathop{\longrightarrow}_{a\to
\infty}}\;
 D_n(\mu(x);b,c,d|q)\; {\mathop{\longrightarrow}_{d\to 0}}\;
 d_n(\mu(x);c,b|q)\; {\mathop{\longrightarrow}_{b= 0}}\;
 V_n^{(c)}(q^{-x};q)
$$
from the four-parameter family of $q$-Racah polynomials, which
occupy the upper level in the Askey-scheme of basic hypergeometric
polynomials (see [19], p. 62), down to the one-parameter set of
Al-Salam--Carlitz II polynomials from the second level in the same
scheme. So, the dual big and dual little $q$-Jacobi polynomials
$D_n(\mu(x);a,b,c|q)$ and $d_n(\mu(x);a,b|q)$ should occupy the
fourth and third level in the Askey-scheme, respectively.

The recurrence relations for the polynomials $D_n(\mu (m)\equiv
D_n(\mu (m);a,b,c|q)$ are obtained from the $q$-difference
equation (6.6). It has the form
$$
(q^{-m}-1)(1-abq^{m+1})D_n(\mu (m))=A_nD_{n+1}(\mu (m))
-(A_n+C_n)D_n(\mu (m))+C_nD_{n-1}(\mu (m)),
$$
where
$$
A_n= (1-aq^{n+1})(1-cq^{n+1}),\ \ \ \ C_n=aq(1-q^n)(b-cq^n).
$$

The relation (8.7) leads to the equality
$$
\sum_{m=0}^\infty\,(-1)^m
\frac{(1-abq^{2m+1})(abq;q)_m}{(1-abq)(q;q)_m}\,
q^{m(m-1)/2}\,D_n(\mu (m);a,b,c|q)\,D_{n'}(\mu (m);
b,a,abq/c|q)=0.    \eqno (8.12)
$$
We give an alternative proof of this result in Appendix.

Note that from the expression (8.8) for the dual big $q$-Jacobi
polynomials $D_n(\mu (m); a,b,c|q)$ it follows that they possess
the symmetry property
$$
D_n(\mu (m); a,b,c|q)=D_n(\mu (m); ab/c,c,b|q). \eqno (8.13)
$$

The set of functions (8.1) and (8.2) form an orthogonal basis in
the Hilbert space ${\frak l}^2$ of functions, defined on the set
of points $m=0,1,2,\cdots$, with the scalar product
$$
\langle f_1,f_2\rangle = \sum_{m=0}^\infty \rho (m) \, f_1(m)\,
\overline{f_2(m)},
$$
where $\rho (m)$ is the same as in formulas (8.5)--(8.7).
Consequent from this fact, one can deduce (in the same way as in
the case of dual little $q$-Jacobi polynomials) that {\it the dual
big $q$-Jacobi polynomials $D_n(\mu (m); a,b,c|q)$ correspond to
indeterminate moment problem and the orthogonality measure for
them, given by formula (8.9), is not extremal}.

It is difficult to find extremal measures. As far as we know,
explicit forms of extremal measures have been constructed only for
the $q$-Hermite polynomials with $q>1$ (see [31]).
\bigskip

\noindent {\bf 9. Generating functions}
\bigskip

Generating functions are known to be of great importance in the
theory of orthogonal polynomials (see, for example, [16]). For the
sake of completeness, we briefly discuss in this section some
instances of linear generating functions for the dual $q$-Jacobi
polynomials $D_n(\mu(x);a,b,c|q)$ and $d_n(\mu(x);a,b|q)$. To
start with, let us consider a generating-function formula
$$
 \sum _{n=0}^\infty \frac{(aq;q)_n}{(q;q)_n}\, t^n D_n(\mu(x);a,b,c|q)
=\frac{(aqt;q)_\infty}{(t;q)_\infty}  \left. {}_2\phi_2 \left(
{q^{-x},\; abq^{x+1}
  \atop abq/c,\;  aqt} \right| q,\; aqt/c   \right) ,  \eqno (9.1)
$$
where $|t|<1$ and, as before, $\mu(x)=q^{-x}+abq^{x+1}$. To verify
(9.1), insert the explicit form (8.8) of the dual big $q$-Jacobi
polynomials
$$
D_n(\mu(x);a,b,c|q)=\sum_{k=0}^n
\frac{(q^{-x},abq^{x+1},q^{-n};q)_k}{(aq,abq/c,q;q)_k} \left(
\frac{aq^{n+1}}{c}\right) ^k
$$
into the left side of (9.1) and interchange the order of
summation. The subsequent use of the relations
$$
(a;q)_{m+k}=(a;q)_m(aq^m;q)_k=(a;q)_k(aq^k;q)_m,
$$  $$
(q^{-m-k};q)_k=(-1)^kq^{-mk-k(k+1)/2}(q^{m+1};q)_k
$$
(see [15], Appendix I) simplifies the inner sum and enables one to
evaluate it by the $q$-binomial formula (4.9). This gives the
quotient of two infinite products in front of ${}_2\phi_2$ on the
right side of (9.1), times $(aqt;q)_k^{-1}$. The remaining sum
over $k$ yields ${}_2\phi_2$ series itself.

As a consistency check, one may also obtain (9.1) directly from
the generating function for the $q$-Racah polynomials
$R_n(\mu(x);\alpha,\beta,\gamma,\delta |q)$ (see formula (3.2.13)
in [19]) by setting $\alpha=q^{-N-1}$ and sending $N\to \infty$.
This results in the relation
$$
 \sum _{n=0}^\infty \frac{(aq;q)_n}{(q;q)_n}\, t^nD_n(\mu(x);a,b,c|q)
=\frac{(aq^{x+1}t;q)_\infty}{(t;q)_\infty}  \left. {}_2\phi_1
\left( {q^{-x},\; c^{-1}q^{-x}
  \atop abq/c} \right| q,\; atq^{x+1}   \right) .  \eqno (9.2)
$$
The left side of (9.2) depends on the variable $x$ by dint of the
combination $\mu(x)=q^{-x}+abq^{x+1}$. Off hand, it is not evident
that the right side of (9.2) is also a function of the lattice
$\mu(x)$. Nevertheless, this is the case. Moreover, the right
sides of (9.1) and (9.2) are equivalent: this fact is known in the
theory of special functions as Jackson's transformation
$$
{}_2\phi_1(a,b;\; c;\; q,z)=\frac{(az;q)_\infty}{(z;q)_\infty}\,
{}_2\phi_2(a,c/b;\; c,az;\; q,bz)
$$
(see, for example, [15]).

The symmetry property (8.13) of the dual big $q$-Jacobi
polynomials $D_n(\mu(x);a,b,c|q)$, combined with (9.1), generates
another relation
$$
 \sum _{n=0}^\infty \frac{(abq/c;q)_n}{(q;q)_n} t^nD_n(\mu(x);a,b,c|q)
=\frac{(abqt/c;q)_\infty}{(t;q)_\infty}  \left. {}_2\phi_2 \left(
{q^{-x},\; abq^{x+1}
  \atop aq,\;  abqt/c} \right| q,\; aqt/c   \right)
$$  $$
=\frac{(abtq^{x+1}/c;q)_\infty}{(t;q)_\infty}  \left. {}_2\phi_1
\left( {q^{-x},\; b^{-1}q^{-x}
  \atop aq} \right| q,\; abtq^{x+1}/c   \right) .  \eqno (9.3)
$$

Similarly, a generating function for the dual little $q$-Jacobi
polynomials has the form
$$
 \sum _{n=0}^\infty \frac{(bq;q)_n}{(q;q)_n} (at)^nd_n(\mu(x);a,b|q)
=\frac{(tq^{-x},abtq^{x+1};q)_\infty}{(at,t;q)_\infty}. \eqno
(9.4)
$$
One can verify (9.4) directly by inserting the explicit form (5.4)
of $d_n(\mu(x);a,b|q)$ into the left  side of (9.4) and repeating
the same steps as in the case of deriving (9.1). This will lead to
the expression
$$
\frac{(abqt;q)_\infty}{(at;q)_\infty}  \, {}_2\phi_1 ( q^{-x},\;
abq^{x+1}; \; abqt; \; q, t )
$$
and it remains only to employ Heine's summation formula (1.5.1)
from [15]. After a simple rescaling of the parameters the
generating function (9.4) coincides with that, obtained earlier in
[24].

The simplest way of obtaining (9.4) is to send $c\to 0$ in both
sides of (9.2): the ${}_2\phi_1$ series on the right side of (9.2)
reduces to ${}_1\phi_0 (q^{-x}; -;\; q,t/a)$, which is evaluated
by the $q$-binomial formula (4.9).

Finally, when the parameter $b$ vanishes, (9.4) reduces to the
known generating function
$$
\sum _{n=0}^\infty
\frac{(-1)^nq^{n(n-1)/2}}{(q;q)_n}t^nV^{(a)}_n(q^{-x};q)
=\frac{(tq^{-x};q)_\infty}{(at,t;q)_\infty}
$$
for the Al-Salam--Carlitz II polynomials (see (3.25.11) in [19]).
\bigskip

\noindent {\bf 10. Concluding remarks}
\bigskip

To summarize, we have attempted to carry over, by using the
discrete series representations of the quantum algebra $U_q({\rm
su}_{1,1})$, an idea of the duality of polynomials, orthogonal on
a finite set of points, to the case of big and little $q$-Jacobi
polynomials.

In fact, we used this idea in [6] and [8] to show that
Al-Salam--Carlitz II polynomials are dual with respect to little
$q$-Laguerre polynomials and $q$-Meixner polynomials are dual to
big $q$-Laguerre polynomials, respectively.

This approach may be further explored and applied to other
$q$-polynomials families. In particular, we know that it is
possible to show, by using certain irreducible representations of
the quantum algebra $U_q({\rm su}_{1,1})$ (which are not
$*$-representations of $U_q({\rm su}_{1,1})$), that polynomials,
dual to $q$-Charlier polynomials (see formula (3.23.1) in [19]),
are Al-Salam--Carlitz polynomials I (see formula (3.24.1) in [19])
and polynomials, dual to the alternative $q$-Charlier polynomials
$$
K_n(x;a;q):= {}_2\phi_1 (q^{-n},-aq^n;\; 0;\; q;qx), \ \ \
n=0,1,2,\cdots ,
$$
(see formula (3.22.1) in [19]), are the polynomials
$$
d_m(\mu(n);a;q)={}_3\phi_0 (q^{-n},-aq^n, q^{-m};\; -;\;
q,-q^m/a), \ \ \ m=0,1,2,\cdots ,
$$
where $\mu (n)=q^{-n}-aq^n$. The orthogonality relation for them
has the form
$$
\sum_{n=0}^\infty \frac{(-a;q)_n(1+aq^{2n})a^n}{(q;q)_n}\,
q^{(3n-1)n/2} d_m(\mu(n)) d_{m'}(\mu(n))=\frac{(-a;q)_\infty
(q;q)_m}{a^mq^{m(m+1)/2}}\, \delta_{mm'},\ \ \ \ a>0.
$$
Proofs of these statements will be given in a separate
publication.
\bigskip

\noindent{\bf Acknowledgments}

\medskip

Discussions with Ch. Berg, J. Christiansen, M. Ismail, T.
Koornwinder, H. Rosengren, and P. Terwilliger are gratefully
acknowledged. This research has been supported in part by the
SEP-CONACYT project 41051-F and the DGAPA-UNAM project IN112300
"Optica Matem\'atica". A.~U.~Klimyk acknowledges the Consejo
Nacional de Ciencia y Technolog\'{\i}a (M\'exico) for a C\'atedra
Patrimonial Nivel II.
\bigskip

\noindent {\bf Appendix}
\bigskip

In this appendix we prove the summation formula
$$
\sum _{n=0}^\infty \frac{(abq,bq;q)_n}{(aq,q;q)_n}
\frac{1-abq^{2n+1}}{1-abq}
a^nq^{n^2}=\frac{(abq^2;q)_\infty}{(aq;q)_\infty}. \eqno (A.1)
$$
First of all, observe that when $b=0$ this relation reduces to
$$
\sum_ {n=0}^\infty \frac{a^nq^{n^2}}{(aq,q;q)_n}=
\frac{1}{(aq;q)_\infty} ,
$$
which is a well-known limiting form of Jacobi's triple product
identity (see [15], formula (1.6.3)).

One can employ an easily verified relation
$$
\frac{(aq,-aq;q)_n}{(a,-a;q)_n}=\frac{1-a^2q^{2n}}{1-a^2}
 \eqno (A.2)
$$
in order to express the infinite sum in (A.1) in terms of a
very-well-poised ${}_4\phi_5$ basic hypergeometric series. This
results in
$$
\sum _{n=0}^\infty \frac{(abq,bq;q)_n}{(aq,q;q)_n}
\frac{1-abq^{2n+1}}{1-abq} a^nq^{n^2} ={}_4\phi_5 \left( \left.
{abq,\; bq,\; q\sqrt{abq},\; -q\sqrt{abq} \atop  aq,\;
\sqrt{abq},\; -\sqrt{abq},\; 0,\; 0} \right|
 q,aq \right) . \eqno (A.3)
$$
The next step is to utilize a limiting case of Jackson's sum of a
terminating very-well-poised balanced ${}_8\phi_7$ series,
$$
{}_6\phi_5 \left( \left. {a,\;  q\sqrt{a},\; -q\sqrt{a},\; b,\;
c,\; d \atop
             \sqrt{a},\;  -\sqrt{a},\; aq/b,\; aq/c,\; aq/d} \right|
 q,\frac{aq}{bcd} \right)
 = \frac{(aq,aq/bc,aq/bd,aq/cd;q)_\infty}
{(aq/b,aq/c,aq/d,aq/bcd;q)_\infty} , \eqno (A.4)
$$
which represents a $q$-analogue of Dougall's formula for a
very-well-poised 2-balanced ${}_7F_6$ series. When the parameters
$c$ and $d$ tend to infinity, from (A.4) it follows that
$$
{}_4\phi_5 \left( \left. {a,\; q\sqrt{a},\; -q\sqrt{a},\; b \atop
    \sqrt{a},\;  -\sqrt{a},\; aq/b,\;  0,\; 0} \right|
 q,\frac{aq}{b} \right)  = \frac{(aq;q)_\infty}{(aq/b;q)_\infty} .
\eqno (A.5)
$$
To verify this, one needs only to use the limit relation
$$
\lim_{c,d\to \infty} (c,d;q)_n\left( \frac{aq}{bcd}\right)^n =
q^{n(n-1)} \left( \frac{aq}{b}\right) ^n .
$$
With the substitutions $a\to abq$ and $b\to bq$ in (A.5), one
recovers the desired identity (A.1).

Similarly, when $d\to \infty$ we derive from (A.4) the identities
$$
\sum_{n=0}^\infty\,\frac{(1-abq^{2n+1})(aq,abq/c,abq;q)_n}
{(1-abq)(bq,cq,q;q)_n (-a/c)^n} q^{n(n-1)/2}=
\frac{(abq^2,c/a;q)_\infty}{(bq,cq;q)_\infty}\, ,    \eqno (A.6)
$$  $$
\sum_{n=0}^\infty\,\frac{(1-abq^{2n+1})(abq,bq,cq;q)_n}
{(1-abq)(aq,abq/c,q;q)_n (-c/a)^n} q^{n(n-1)/2}
=\frac{(abq^2,a/c;q)_\infty}{(aq,abq/c;q)_\infty}\, .  \eqno (A.7)
$$
They have been used in section 7.

We conclude this appendix with the following remark. There is
another proof of the identity (8.12), based on vital use of the
same summation formula (A.4). Actually, a relation may be derived,
which is somewhat more general than (8.12). Indeed, consider the
function
$$
\eta_k(a;q):=\sum_{n=0}^\infty (-1)^nq^{n(n-1)/2}
\frac{1-aq^{2n+1}}{1-aq} \frac{(aq;q)_n}{(q;q)_n} \mu^k(n;a)
 \eqno (A.8)
$$
for arbitrary nonnegative integers $k$, where the $q$-quadratic
lattice $\mu(n; a)$ is defined as before:
$$
\mu(n;a):= q^{-n}+aq^{n+1}.  \eqno (A.9)
$$
We argue that all $\eta_k(a;q)=0$, $k=0,1,2,\cdots$. To verify
that, begin with the case when $k=0$ and employ relation (A.2) to
show that
$$
\eta_0(a;q)=\left. {}_3\phi_3 \left( {q\sqrt{aq},\; -q\sqrt{aq},\;
aq     \atop \sqrt{aq},\;  -\sqrt{aq},\; 0} \right| q, 1 \right) .
$$
The summation formula (A.4) in the limit as $d\to \infty$ takes
the form
$$
\left. {}_5\phi_5 \left( {a,\; q\sqrt{a},\; -q\sqrt{a},\; b,\; c
  \atop \sqrt{a},\;  -\sqrt{a},\; aq/b,\; aq/c,\; 0} \right| q,\;
  \frac{aq}{bc}
  \right) =\frac{(aq,aq/bc;q)_\infty}{(aq/b,aq/c;q)_\infty}\ .
 \eqno (A.10)
$$
In the particular case when $bc=aq$ this sum reduces to
$$
\left. {}_3\phi_3 \left( {a,\; q\sqrt{a},\; -q\sqrt{a}
  \atop 0,\; \sqrt{a},\;  -\sqrt{a}} \right| q, 1
  \right) =\frac{(aq,1;q)_\infty}{(b,c;q)_\infty}=0 \,,
$$
since $(z;q)_\infty =0$ for $z=1$. Consequently, the function
$\eta_0(a;q)$ does vanish.

For $k=1,2,3,\cdots$, one can proceed inductively. Employ the
relation $q\mu (n+1;a)=\mu (n;q^2a)$ to show that
$$
\eta_{k+1}(a;q)=(1+aq)\eta_k(a;q)-q^{-k-1}(1-aq^2)(1-aq^3)
\eta_k(aq^2;q).
$$
So, one obtains that indeed all $\eta_k(a;q)$, $k=0,1,2,\cdots$,
vanish. The identity (8.12) is now an easy consequence of this
statement if one takes into account that a product of the two
polynomials $D_n(\mu(m); a,b,c|q)$ and $D_{n'}(\mu(m);
b,a,abq/c|q)$ in (8.12) is some polynomial in $\mu(m)$ of degree
$n+n'$. This completes the proof of (8.12), which is independent
of the one, given in section 8.

\end{document}